\documentclass[12pt]{article}

\usepackage{amsmath,amssymb}

\newtheorem{THM}{Theorem}[section]

\begin{document}

\begin{center}
{\large\bf On meromorphic functions without JULIA directions}\\
[.4cm]
{Peter Tien-Yu  Chern } \footnote{
\indent This paper was supported in part by NSC (R.O.C)under the
grant NSC 95-2115-M-214-003 and was addressed in ICM August 2006,
Madrid, Spain}
\end{center}

\bigskip

\begin{abstract}
It is proved that for any positive number $\lambda$, $1<\lambda<2$;
there exists a meromorphic function $f$ with logarithmic order
$\lambda$= $\displaystyle\limsup_{r\to+\infty}\frac{\log
T(r,f)}{\log\log r}$  such that $f$ has no Julia directions, where
$T(r,f)$ is the Nevanlinna characteristic function of $f$. (Note
that A. Ostrowski has proved a like-result for $\lambda=2$ in 1926.)
\end{abstract}

\section{Introduction}

In 1919 G. Julia [J] essentially obtain the following: Let $f(z)$ be
a function in the complex plane $\mathbb C$. Let $z_0$ be an
essential singularity of $f(z)$. Then there exists at least a ray
$\Delta (\theta)=\{z_0+re^{i\theta}:\ r\geq 0\}$, where
$0\leq\theta<2\pi$, emanating from $z_0$, such that in every angular
region $\theta-\varepsilon<\arg(z-z_0)<\theta+\varepsilon$,
$\varepsilon>0$, $f(z)$ assumes, infinitely often, every extended
complex value with at most two exceptions.

The ray above is called {\it a Julia direction for $f(z)$.} Every
transcendental entire function $f$ has $\infty$ as an essential
singularity, hence by above Julia's result, we deduce that $f(z)$
has a Julia direction.

Let $\phi(r)$ be a positive increasing function defined in
$(0,+\infty)$. The logarithmic order of $\phi(r)$ is the number
\begin{align*}
\limsup_{r\to +\infty}\frac{\log\phi(r)}{\log\log(r)}.
\end{align*}
Let $f$ be a meromorphic function in the complex plane, we say $f$
has logarithmic order $\lambda$ if $T(r,f)$ has logarithmic order
$\lambda$, where $T(r,f)$ is the Nevanlinna characteristic function of f.

The aim of this article is to prove the following

\begin{THM}
For each $\lambda$, $1<\lambda<2$, there exists a transcendental
meromorphic function $f$ of logarithmic order $\lambda$, such that
$f$ has no Julia direction.
\end{THM}

\section{Results on Linear Fractional Transformations}

For any finite complex number $\alpha$, we define a linear
fractional transformation $\omega_\alpha:\mathbb C\to\mathbb C,$ by
$\omega_\alpha(z)=\frac{\alpha+z}{\alpha-z}$.

In case $\alpha$ is positive, we obviously have the following
properties:\\
\indent\indent If $Re\ z>0$, then $|\omega_\alpha(z)|>1$;\\
\indent\indent if $Re\ z<0$, then $|\omega_\alpha(z)|<1$,\\
and
\indent if $Re\ z=0$, then $|\omega_\alpha(z)|=1$.

For $K$,$0<K<1$, we put
\begin{align*}
F_\alpha(K)=\left\{z:\ z\in\mathbb C,\ |\omega_\alpha(z)|<K\right\}.
\end{align*}

Then $F_\alpha(K)$ is a disk which is symmetric with respect to the
$x$-axis in the left half plane $\{z:\ Re\ z<0\}$ with center
$(O_\alpha,0)$ where $O_\alpha=-\frac{K^2+1}{1-K^2}\alpha$ and
radius $R_\alpha=\frac{2K}{1-K^2}\alpha$. The near point of
$F_\alpha(K)$  in the $x$-axis is $(X_\alpha(K),0)$,
where $X_\alpha(K)=O_\alpha+R_\alpha=-\frac{1-K}{1+K}\alpha$,
and the farthest point of $E_\alpha(K)$ in the $x$-axis is
$(\chi_\alpha(K),0)$,
where $\chi_\alpha(K)=O_\alpha-R_\alpha=-\frac{1+K}{1-K}\alpha$.
 Since $|O_\alpha|>R_\alpha$, so $F_\alpha(K)\varsubsetneq\{z:\ |\arg
z-\pi|<\frac{\pi}{4}\}$.

\section{The proof of the Theorem}

Let $A_1=1$, $A_n=\exp\left(n^{\frac{1}{\lambda-1}}\right)$ for
$n=2,3,4,\cdots$. For $n=1,2,\cdots$, we set
\begin{align*}
E_n=E_{A_n}=\left\{z:\ z\in\mathbb C,\ |\omega_{A_n}(z)|<\frac
13\right\},
\end{align*}
and
\begin{align*}
F_n=F_{A_n}=\left\{z:\ z\in\mathbb C,\
|\omega_{A_n}(z)|<I_n\right\},
\end{align*}
where $\displaystyle I_n=\frac{1-\frac{1}{n+1}}{1-\frac
1n}=\frac{n^2+2n}{(n+1)^2}$. we have $E_{A_n}\subset\{z:\ |\arg
z-\pi|<\frac{\pi}{4}\}$, and $E_{A_n}\varsubsetneq F_{A_n}$. Further
we have
\begin{align*}
X_{A_n}(I_n)=\frac{-A_n}{2n^2+4n+1}
\end{align*}
and
\begin{align*}
\chi_{A_n}(I_n)=-(2n^2+4n+1)A_n.
\end{align*}

Applying L H\^opital's rule, we can obtain that there is a positive
integer $n_0$ depending only on $\lambda$, such that if $n>n_0$, then
\begin{align*}
\chi_{A_{n+1}}(I_{n+1})<X_{A_{n+1}}(I_{n+1})<\chi_{A_n}(I_n)<X_{A_n}(I_n),
\end{align*}
this is correct, since $1<\lambda<2$. Hence $\mathfrak
F=\left\{F_{A_n}(I_n)\right\}_{n=n_0+1}^{+\infty}$ is a family of
mutually disjoint open sets in $\mathbb C$.

Now we define a function $f:\mathbb C\to \mathbb C$ by
\begin{align*}
f(z)=\prod_{j=n_0+1}^{+\infty}\omega_{A_j}(z),\ \ \
A_j=\exp\left(j^{\frac{1}{\lambda-1}}\right).
\end{align*}
We shall show that $T(r,f)$ has logarithmic order $\lambda$ and $f$
has no Julia direction. We first show that $f$ has no Julia
direction.

For each $\displaystyle
z\in\left(\bigcup_{n=1+n_0}^{+\infty}E_{A_n}\right)^c=\bigcap_{n=1}^{+\infty}E_{A_{n_0+n}}^c$,
\begin{align*}
|f(z)|\geq \frac
13\times\prod_{n=1+n_0}^{+\infty}\frac{1-\frac{1}{n+1}}{1-\frac
1n}=\frac{1/3}{1+n_0}.
\end{align*}

Since $\mathfrak F$ is a family of disjoint sets, for this $z$ there
exists at most one positive integer $k_0$ such that $z\in
F_{A_{k_0}}$.

In
$\displaystyle\left(\bigcup_{n=1+n_0}^{+\infty}E_{A_n}\right)^c$,
$f$ omits values in the disk
$D\left(0,\frac{1/3}{n_0+1}\right)=\{z:\ |z|<\frac{1/3}{n_0+1}\}$.
Since $E_{A_n}\subset\left\{z:\ |\arg
z-\pi|<\frac{\pi}{4}\right\}$. $f$ has no Julia direction in
$\Omega_{1}=\displaystyle\{z:\ |\arg z|< (3/4)\pi\} $

For $Re\ z=x<0$, $|\omega_{A_n}(z)|<1$, it follows that
\begin{align*}
|f(z)|<1\ \ \ \mbox{for}\ z\in\Omega_2=\left\{z:\ |\arg
z-\pi|<\frac{\pi}{2}\right\}.
\end{align*}
This implies that $f$ has no Julia direction in  $\Omega_2$.
Since $\mathbb C$ is the union of $\Omega_1$ and $\Omega_2$,
so f has no Julia direction in the complex plane $\mathbb C$.
\

We next shall prove that $T(r,f)$ has logarithmic order $\lambda$. Since the
sets $\{A_n\}_{n=1}^{+\infty}$ and $\{-A_n\}_{n=1}^{+\infty}$ both
have logarithmic convergence exponent $\lambda-1$, it follows from
Theorem 3.1 and Theorem 4.1 of [C2] that $N(r,f=0)$ and
$N(r,f=\infty)$ both have logarithmic order $\lambda$. Since $f$ is
a transcendental function, so the lower logarithmic order, say,

\begin{align*}
\mu=\liminf_{r\to+\infty}\frac{\log T(r,f)}{\log\log r}\geq 1,
\end{align*}

and $\lambda<2$, we have $\lambda-\mu<1$, it follows from Theorem
7.1 of [C2],
\begin{align*}
T(r,f)\leq N(r,f=0)+N(r,f=\infty)+o(T(r,f)),
\end{align*}
further, $N(r,f=0)$ and $N(r,f=\infty)$ both have logarithmic order
$\lambda$, hence $T(r,f)$ has logarithmic order $\lambda$. This
completes the proof of the Theorem.

{\bf Remark. 1.} For $\lambda>2$, any meromorphic function with
logarithmic order $\lambda$ must have a Julia direction, since in
[C1] it is proved that if f is a meromorphic function with finite
logarithmic order, and satisfies the growth condition
\begin{align*}
\limsup_{r\to+\infty}\frac{T(r,f)}{(\log r)^2}=+\infty,
\end{align*}
then $f$ has a Borel direction with logarithmic order $\lambda-1$,
and a Borel direction of $f$ must be a Julia direction of $f$. For
$\lambda=2$, an example has been given by A.Ostrowski in 1926[O].

Before the end of this article, I submit two research problems.

\indent 1. Is there a transcendental meromorphic function $f$ with
logarithmic order $1$, such that $f$ has no Julia direction.

\indent 2. For any positive number $\lambda$, with $1<\lambda<2$,
does there exist a function $f$ with logarithmic order $\lambda$ and
a ray $\Delta(\theta)=\{re^{i\theta}|r\geq 0\}$, such that
$\Delta(\theta)$ is a Julia direction of $f$, but $\Delta(\theta)$
is not a Borel direction for $f$?

\bigskip

Department of Applied Mathematics

I-Shou University

Ta-Hsu,Kaohsiung 840 Taiwan R.O.C.

E-Mail: tychern@mail.isu.edu.tw

\end{document}